\documentclass{agtart_a}
\pdfoutput=1

\usepackage{amscd}
\usepackage[all]{xy}

\title{Cosmetic surgeries on genus one knots}

\author{Jiajun Wang}
\givenname{Jiajun}
\surname{Wang}
\address{Department of Mathematics\\
University of California at Berkeley\\\newline
Berkeley, CA 94720\\USA}
\email{wang@math.berkeley.edu}
\urladdr{}

\volumenumber{6}
\issuenumber{}
\publicationyear{2006}
\papernumber{55}
\startpage{1491}
\endpage{1517}

\doi{}
\MR{}
\Zbl{}


\received{26 February 2006}
\revised{26 July 2006}
\accepted{17 August 2006}
\published{4 October 2006}
\publishedonline{4 October 2006}
\proposed{}
\seconded{}
\corresponding{}
\editor{}
\version{}

\arxivreference{math.GT/0512253}

\let\xysavmatrix\xymatrix
\def\xymatrix{\disablesubscriptcorrection\xysavmatrix}
\AtBeginDocument{\let\bar\wbar\let\tilde\wtilde\let\hat\what}
\makeop{rk}
\newcommand{\abs}[1]{\lvert#1\rvert}
\newcommand{\babs}[2]{#1\lvert#2#1\rvert}

\makeatletter
\def\cnewtheorem#1[#2]#3{\newtheorem{#1}{#3}[section]
\expandafter\let\csname c@#1\endcsname\c@thm}
\makeatother

\newtheorem{thm}{Theorem}[section]
\cnewtheorem{lemma}[thm]{Lemma}
\cnewtheorem{prop}[thm]{Proposition}
\cnewtheorem{cor}[thm]{Corollary}
\cnewtheorem{conj}[thm]{Conjecture}

\newcommand{\bba}{\mathbb{A}}
\newcommand{\bbb}{\mathbb{B}}
\newcommand{\bbc}{\mathbb{C}}
\newcommand{\bbq}{\mathbb{Q}}
\newcommand{\bbx}{\mathbb{X}}
\newcommand{\bbz}{\mathbb{Z}}
\newcommand{\calf}{\mathcal{F}}
\newcommand{\calt}{\mathcal{T}}
\newcommand{\fraks}{\mathfrak{s}}
\newcommand{\frakt}{{\mathfrak{t}}}

\newcommand{\dps}[1]{\displaystyle{#1}}

\newcommand{\spinc}{${\rm Spin}^c$ }
\newcommand{\ch}{{\rm Char}}

\begin{document}

\begin{asciiabstract}
In this paper, we prove that there are no truly cosmetic surgeries
on genus one classical knots. If the two surgery slopes have the same sign,
we give the only possibilities of reflectively cosmetic surgeries. The result
is an application of Heegaard Floer theory and number theory.
\end{asciiabstract}

\begin{abstract}   
In this paper, we prove that there are no truly cosmetic surgeries on genus one classical knots. If the two surgery slopes have the same sign, we give the only possibilities of reflectively cosmetic surgeries. The result is an application of Heegaard Floer theory and number theory.
\end{abstract}

\maketitle

\section{Introduction}

Let $Y$ be a closed oriented three-manifold and $K$ be a framed knot in $Y$.
For a rational number $r$, let $Y_r(K)$ be the resulting manifold
obtained by Dehn surgery along $K$ with slope $r$ with respect to the given framing. For a knot $K$ in $S^3$, we will always use the Seifert framing. Two surgeries along $K$ with distinct slopes $r$ and $s$ are called {\em cosmetic\/} if $Y_r(K)$ and $Y_s(K)$ are homeomorphic. The two surgeries are {\em truly cosmetic\/} if the homeomorphism is orientation preserving and {\em reflectively cosmetic\/} if the homeomorphism is orientation reversing. Cosmetic surgeries $Y_r(K)$ and $Y_s(K)$ are {\em mundane\/} if there is a homeomorphism of the knot exterior taking $r$ to $s$ and {\em exotic\/} otherwise.

Surprisingly, (reflectively) cosmetic surgeries are not rare at all. For an
amphichiral knot $K$ in $S^3$, we always have $S^3_r(K) \cong -S^3_{-r}(K)$.
Interesting examples are given by Mathieu \cite{Mathieu}. He constructed a
series of cosmetic surgeries on the trefoil $T$, namely, for any nonnegative
integer $k$, we have $$S^3_{(18k+9)/(3k+1)}(T)\cong -S^3_{(18k+9)/(3k+2)}(T).$$
The resulting manifold $M_k$ is the Seifert fibred space which has Seifert invariant
$S(0; k{-}3/2; (2,1), (3,1), (3,2))$ in the notation of Jankins and Neumann \cite{Neuman}. By results
in Seifert fibred spaces, $M_k$ admits no orientation-reversing homeomorphisms.
So they are not truly cosmetic. As pointed out in \cite{BleilerHodgsonWeeks},
the existence of Seifert fibred cosmetic surgeries highly depends on the
existence of an exceptional fibre of index 2, which could be recovered by
Rolfsen twist after orientation reversal.

An even more striking example is given by Bleiler, Hodgson and Weeks
\cite{BleilerHodgsonWeeks}. They give an oriented $1$--cusped hyperbolic $3$--manifold
$X$ with a pair of slopes $r_1$ and $r_2$ such that the two surgered manifolds
$X(r_1)$ and $X(r_2)$ are oppositely oriented copies of the lens space
$L(49,18)$ (and there is no homeomorphism $h$ of $X$ such that $h(r_1)=r_2$).
However, it is still unknown whether there are hyperbolic manifolds that could
be obtained via exotic (truly or reflectively) cosmetic surgeries. A result in
\cite{BleilerHodgsonWeeks} indicates that such examples would be unexpected.

There are truly cosmetic surgery on the unknot: $S^3_{\smash{p/q}}(U)\cong
S^3_{\smash{p/p+q}}(U)$ (mundane) and $S^3_{\smash{p/q_1}}(U)\cong S^3_{\smash{p/q_2}}(U)$ (exotic)
where $q_1q_1\equiv1\pmod{p}$. However, so far, there are no known truly
cosmetic surgeries on a nontrivial knot. In fact, recall the following
conjecture, Problem 1.81(A) in Kirby's problem list:

\begin{conj}[Cosmetic Surgery Conjecture \cite{Kirby}] \label{cosmeticconj}
There are no truly cosmetic surgeries on a knot $K$ in a closed three-manifold
$Y$ provided that the knot complement is not homeomorphic to the solid torus. In
particular, for any classical nontrivial knot $K$ in $S^3$, $S^3_r(K)\cong
S^3_s(K)$ if and only if $r=s.$ \end{conj}

One approach to this conjecture is to get restrictions on manifolds that could
arise as the resulting manifold of cosmetic surgeries. For example, the
three-sphere (Gordon and Luecke \cite{GordonLuecke}) and the manifold
$S^2\times S^1$ (Gabai \cite{Gabai}) could not be obtained via cosmetic
surgeries (truly or reflectively). In \cite{Rong}, Rong classified all Seifert
fibered spaces that could be obtained by cosmetic surgeries on a knot whose
complement is Seifert fibered and boundary incompressible. One remarkable
consequence of Rong's classification is that if a knot $K$ with Seifert fibered
complement admits a pair of cosmetic surgeries yielding a Seifert fibered space,
then it admits infinitely many such cosmetic surgeries. All these cosmetic
surgeries are reflectively cosmetic.

Another approach is to get restrictions on knots in manifolds that admit
cosmetic surgeries. In \cite{Lackenby}, Lackenby showed that for two surgered
manifolds $Y_{\smash{p/q}}(K_1)$ and $Y^\prime_{\smash{p^\prime /q^\prime}}(K_2)$, if $Y$ and
$Y^\prime$ are distinct, or $K_1$ and $K_2$ are distinct, or $p/q\neq
p^\prime/q^\prime$, then we have $Y_{\smash{p/q}}(K_1)\not\cong Y^\prime_{\smash{p^\prime
/q^\prime}}(K_2)$ provided that $\abs{q^\prime}$ is sufficiently large. This
implies that truly cosmetic surgeries are really rare.

Heegaard Floer theory was introduced by Ozsv{\'a}th and Szab{\'o} in
\cite{OSzClosed3M} and other sequels; see \cite{OSzSurvey} for an overview.
Heegaard Floer theory assigns invariants to closed three-manifolds (called
Heegaard Floer homologies) and four-manifolds. And there are also numerous
other invariants, such as knot Floer homologies and Ozsv{\'a}th--Szab{\'o}
contact invariants. Heegaard Floer theory (especially the rational surgery
formula \cite{OSzRationalSurgeries}) gives strong restrictions on (classical)
knots that admit cosmetic surgeries. Roughly speaking, an {\em L-space\/} is a
closed $3$--manifold having the same Heegaard Floer homology as some lens space.
In particular, we have the following theorem.

\begin{thm}[Ozsv{\'a}th--Szab{\'o} \cite{OSzRationalSurgeries}]
\label{hfkofcosmeticknot} Let $K\subset S^3$ be a knot. If $r, s$ are distinct
rational numbers and $S^3_r(K)\cong\pm S^3_s(K)$, then either $S^3_r(K)$ is an
$L$--space or $r$ and $s$ have opposite signs. If $K$ is of Seifert genus one,
and $S_r^3(K)\cong S^3_s(K)$ with $r\neq s$, then $S^3_r(K)$ is an $L$--space and
$K$ has the same knot Floer homology (and hence the same Alexander polynomial)
as the trefoil knot $T$.
\end{thm}

In this paper, based on the above theorem, we prove \fullref{cosmeticconj} for genus one classical knots. Specifically, our main result
is the following theorem.

\begin{thm}\label{genusonecosmetic} For two Dehn surgeries with slopes $p/q$ and
$p/q^\prime$ on a Seifert genus one knot $K$ in $S^3$, let $S^3_{\smash{p/q}}(K)$ and
$S^3_{\smash{p/q^\prime}}$ be the resulting manifolds. Then
\begin{enumerate}
\item $S^3_{\smash{p/q}}(K)\cong S^3_{\smash{p/q^\prime}}(K)$ if and only if $q=q^\prime$.
\item If $S^3_{\smash{p/q}}(K)\cong-S^3_{\smash{p/q^\prime}}(K)$ and the slopes $p/q$ and $p/q^\prime$ have the
same sign, then either $K$ has the same knot Floer homology as the right-handed
trefoil knot and $p=18k{+}9$, $\{q, q^\prime\}= \{3k{+}1, 3k{+}2\}$ for some
nonnegative integer $k$, or  $K$ has the same knot Floer homology as the left-handed trefoil knot and $p=18k{+}9$, $\{q, q^\prime\}= \{-(3k{+}1), -(3k{+}2)\}$ for
some nonnegative integer $k$.
\end{enumerate}
\end{thm}

Our theorem guarantees that Mathieu's examples are the only cosmetic surgeries
on the right-handed trefoil. In fact in their survey paper \cite{OSzSurvey},
Ozsv{\'a}th and Szab{\'o} implicitly conjectured that knot Floer homologies
detect fibered of knots in $S^3$, namely, a knot $K$ in $S^3$ is fibred if
and only its knot Floer homology in the topmost nontrivial filtration is
isomorphic in $\bbz$. Several classes of knots have been verified for this
conjecture: pretzel knots in Ozsv{\'a}th and Szab{\'o} \cite{OSzMutation}, some cable knots in
Hedden \cite{Hedden}, Whitehead double knots in Eftekhary \cite{Eftekhary}, closed $3$--braids in
Ni \cite{Ni3Braids}, etc. See Ni \cite{NiSutured} for  theoretical evidence. So
by \fullref{hfkofcosmeticknot}, we expect that the reflectively cosmetic
surgeries on the trefoil knot are the only cosmetic surgeries on genus one
knots if the slopes have the same sign. (Of course on amphichiral knots, such
as the figure eight knot, there are many reflectively cosmetic surgeries with
slopes in opposite signs by letting $r=-s$.)

The ingredients in our proof of \fullref{genusonecosmetic} are Heegaard
Floer theory, Reidemeister torsion, the Casson--Walker invariant and number theory
on cyclotomic numbers. Truly cosmetic surgeries of slopes with opposite signs
on genus one knots are ruled out by \fullref{hfkofcosmeticknot}. For
cosmetic surgeries (truly or reflectively) with slopes of the same sign,
\fullref{hfkofcosmeticknot} says that such genus one knots must have the
same knot Floer homology as the trefoil knot. Then we use Reidemeister torsion
together with number theory to get restrictions on the slopes that might appear
as cosmetic surgeries. Finally, we use Heegaard Floer theory and the Casson--Walker
invariant to rule out the remaining possibilities.

The paper is organized as follows. In \fullref{hfhprelim}, we recall some
preliminaries on Heegaard Floer homology. In \fullref{cyclotomic}, we
establish a number theoretic result regarding cyclotomic numbers which is useful
in \fullref{cosmetic}. In \fullref{cosmetic}, we use number theory and
Heegaard Floer homology to prove \fullref{genusonecosmetic}.

{\bf Remark}\qua After the submission of this paper, P\,Ghiggini \cite{Ghiggini} proved the fibred knot conjecture for genus one knots and Y\,Ni \cite{Nifibred} proved the general case. Hence our theorem could be obtained using their results without number theory.

{\bf Acknowledgements}\qua The paper was completed while the author was an
exchange student at Columbia University. The author is grateful to Robion Kirby
and Peter Ozsv\'ath for their support and many
helpful discussions. We thank Elisenda Grigsby, Ming-Lun Hsieh, Walter Neumann, Bjorn Poonen, Eric Urban and Shaffiq Welji for interesting discussions. Also, we would like to thank the referee for suggestions and comments on the first draft.

\section{Preliminaries in Heegaard Floer theory}\label{hfhprelim}

We review some of the materials and notation in Heegaard Floer theory which
will be used in this article.

\subsection[Correction terms and spin-c structures]{Correction terms and \spinc structures}

For a rational homology three-sphere $Y$ with \spinc structure $\fraks$, the
minimal absolute $\bbq$ grading of torsion-free elements in the image of
$HF^\infty(Y,\fraks)$ in $HF^+(Y, \fraks)$ is an invariant of $(Y,\fraks)$,
called the {\em correction term\/} of $(Y,\fraks)$ and denoted by $d(Y,\fraks)$.
See Ozsv{\'a}th and Szab{\'o} \cite{OSzAbsoluteGrade} for details. One basic property of the correction
terms is
\begin{equation}\label{CorrectionTermProperties}
d(Y,\fraks)=d(Y,\bar\fraks),\quad d(Y,\fraks)=-d(-Y,\fraks).
\end{equation}
The correct term contains information of the intersection forms of
four-manifolds which bound $Y$. Specifically, if $X$ is an oriented
four-manifold whose intersection form is negative-definite, and $\frakt$ is any
\spinc structure on $X$ whose restriction to $Y=\partial X$ is $\fraks$, then we
have the inequality
\begin{equation}
\frac{c_1(\frakt)^2+\rk H^2(X;\bbz)}4\leq d(Y,\fraks).
\end{equation}
In \cite{OSzPlumb}, Ozsv{\'a}th and Szab{\'o} showed that if $X$ is a plumbed
four-manifold with a negative-definite graph $G$ having at most two ``bad
vertices," then the correction term is a sharp bound in the above inequality.
(A vertex $v$ is {\em bad\/} if the number of edges containing $v$ is bigger than
minus the Euler number of the corresponding disk bundle.) Specifically, for each
\spinc structure $\fraks$ of $Y$, the correction term $d(Y,\fraks)$ satisfies
$$d(Y,\fraks)=\max_{K\in{\rm Char}(X, \fraks)}\frac{K^2+\abs{G}}4,$$
where $\abs{G}$ is the number of vertices of $G$, and $\ch(X, \fraks)$ is the
set of first Chern classes of those \spinc structures of $X$ whose restriction
to $Y$ is $\fraks$. This gives a practical way to compute the correction terms
for a three-manifold $Y$ which could be realized as the boundary of such a
plumbed four-manifold.

In particular, a linear plumbing with all self-intersections at most $-2$ has no bad
vertices. In what follows, we just consider such linearly plumbed
four-manifolds. A linearly plumbed four-manifold could be denoted by $(a_0, a_1,
\cdots, a_n)$ ($a_i\leq-2$), where $[S_i]\cdot[S_i]=a_i$ and
$[S_i]\cdot[S_{i+1}]=1$ (we denote by $[S_i]$ the homology class of the $i$--th
exceptional sphere). Let $\phi_i$ denotes the Poincar{\'e} dual of $S_i$.
Note that if $(a_0, a_1, \cdots, a_n)$ gives the continued fraction expansion of
$p/q$, ie,$$
\displaystyle{\frac pq=[a_0,a_1,\cdots,a_n]=a_1-\frac{1}{
a_2-\displaystyle{\frac{1}{\cdots-\displaystyle{\frac{1}{a_n}}}}}},
$$
then the boundary of the corresponding linearly plumbed four-manifold is the
lens space $-L(p,q)=S^3_{\smash{p/q}}(U)$.

Let $X$ be linearly plumbed four-manifold $(a_0,a_1,\cdots,a_n)$ with
$a_i\leq-2$ and $Y$ be the boundary three-manifold of $X$. $Y$ is a rational
homology sphere. The set of first Chern classes of \spinc structures on $X$ is
exactly the set of characteristic classes of $X$. Let $\fraks$ be a \spinc
structure on $Y$ and $\frakt$ be a \spinc structure on $X$ whose restriction on
$Y$ is $\fraks$, then
\begin{equation*}\ch(X,\fraks)=\{c_1(\frakt)+C_0\phi_0+\cdots+C_n\phi_n\ \big|\
\text{$C_i$'s are even integers}\}.
\end{equation*}
If the first homology group of $Y$ is cyclic of an odd order, there is a unique
spin \spinc structure on $Y$, ie, the one whose first Chern class is zero. We
denote this spin \spinc structure by $\fraks_0$. The set ${\rm Char}(X, \fraks_0)$
contains those characteristic classes which could be written as a linear
combination of the $\phi_i$'s.

\subsection{Rational surgery formulas} Fix a nullhomologous knot $K$ in
a closed three-manifold $Y$. Knot Floer homology associates to $K$ a
$\bbz\oplus\bbz$--filtered chain complex $C=CFK^\infty(Y,K)$ which is a
$\bbz[U]$--module generated over $\bbz$ by a set $X$. The filtration function
$\calf\co X\rightarrow\bbz\oplus \bbz$ satisfies that, if $\calf({\bf x})=(i,j)$,
then $\calf(U\cdot{\bf x})=(i-1,j-1)$ and $\calf ({\bf y})\leq\calf({\bf x})$
for all $\bf y$ having nonzero coefficient in $\partial{\bf x}$.

For a subset $S$ of $\bbz\oplus\bbz$ with the property $(i,j)\in S$
implies $(i+1,j),(i,j+1)\in S$, let $C\{S\}$ be the quotient complex of $C$ by
the subcomplex generated by $\bf x$ with $\calf({\bf x})\notin S$. For an
integer $s$, let $A^+_{\smash{s}}(K)=C_K\{\max(i,j-s)\geq0\}$ and $B^+(K)=C\{i\geq0\}$.
There is a canonical chain homotopy equivalence (up to sign) between $B^+$ and
$C\{j\geq0\}$, denoted by $T$. Also, there are two chain maps $v_s^+,
h_s^+\co A^+_s\rightarrow B^+$. The map $v_s^+$ is the projection of $A^+_s$ onto
$B^+$, and $h_s^+$ is the composition map:
$$\begin{CD} A^+_s
@>{\rm projection}>> C\{j\geq s\}  @>{U^s\cdot}>> C\{j\geq0\}@>{T}>>B^+
\end{CD}$$
We will simply write as $v^+$ and $h^+$ when there is no confusion.


Now fix a surgery slope $p/q$ and we suppose $q>0$. Consider the two chain
complexes
$$\bba^+=\bigoplus_{t\in\bbz}(t,A^+_{\lfloor\unfrac{t}{q}\rfloor}),\quad
\bbb^+=\bigoplus_{t\in\bbz}(t,B^+),$$ where $\lfloor x\rfloor$ is the greatest
integer not bigger than $x$. An element of $\bba^+$ could be written as
$\{(t,a_t)\}_{t\in\bbz}$ with $a_t\in A^+_{\lfloor\unfrac tq\rfloor}$. Define a
chain map $D^+_{\smash{p/q}}\co \bba^+\rightarrow\bbb^+$ by
$$D^+_{\smash{p/q}}\{(t,a_t)\}=\{(t,b_t)\},\quad
b_t=v^+(a_t)+h^+(a_{t-p}).$$ Let $\bbx_{\smash{p/q}}^+(K)$ be the mapping cone of
$D^+_{\smash{p/q}}$. Note that since both $A^+_s$ and $B^+_s$ are relatively
$\bbz$--graded, and the two chain maps $h^+$ and $v^+$ respect this relative
grading, the chain complex $\bbx_{\smash{p/q}}^+$ can be relatively $\bbz$--graded with
the convention that the differential $D^+_{\smash{p/q}}$ lowers the grading by one.

$\bbx^+_{\smash{p/q}}$ naturally splits into the direction sum of $p$ subcomplexes
$$\bbx^+_{p,q}=\bigoplus_{i=0}^{p-1}\bbx^+_{i,p/q},$$
where $\bbx^+_{i,p/q}$ is the subcomplex of $\bbx^+_{\smash{p/q}}$ containing all
$A^+_t$ and $B^+_t$ with $t\equiv i\pmod{p}$.

Ozsv{\'a}th and Szab{\'o} showed that the Heegaard Floer homology of a $p/q$
surgered manifold is determined by the mapping cone $\bbx^+_{\smash{p/q}}$.

\begin{thm}[Ozsv{\'a}th and Szab{\'o} \cite{OSzRationalSurgeries}]\label{rationalsurgeryformula}
 Let $K\in
Y$ be a nullhomologous knot and $p,q$ be two coprime integers.
Then there are relatively graded isomorphisms
$$HF^+(Y_{\smash{p/q}}(K))\cong H_*(\bbx^+_{\smash{p/q}}),\quad HF^+(Y_{\smash{p/q}}(K),\fraks_i)
\cong H_*(\bbx^+_{i,p/q}),$$ where $\fraks_i$ is the \spinc structure
corresponding to $i\in\bbz/p\bbz$.
\end{thm}

In fact, the absolute grading on $HF^+(Y_{\smash{p/q}}(K))$ is also determined. The
absolute grading on $\bbx^+_{\smash{p/q}}$ is specified by fixing an absolute grading
on $B^{\smash{+}}_{\lfloor\unfrac tq\rfloor}$ for any $t$. The absolute grading on $B^{\smash{+}}_{
\smash{\lfloor\unfrac tq\rfloor}}$ is determined by a grading on its homology
$H_*(B^{\smash{+}}_{\smash{\lfloor\unfrac tq\rfloor}})\cong\calt^+$, where $\calt^+$ is the
$\bbz[U]$--module $\bbz[U,U^{-1}]/(U\cdot\bbz[U])$. The absolute grading on
$H_*(B^+_{\smash{\lfloor\unfrac tq\rfloor}})\cong\calt^+\cong H_*(\bbx^+_{\smash{i,p/q}})$ is
determined so that its bottom-most nontrivial element is supported in dimension
$d(S^3_{\smash{p/q}}(U),i)=d(-L(p,q),i)$. These correction terms $d(-L(p,q),i)$
are computed by the recursive formula
\cite[Proposition 4.8]{OSzAbsoluteGrade}
\begin{equation}
d(-L(p,q),i)=\left(\frac14-\frac{(2i+1-p-q)^2}{4pq}\right)-d(-L(q,r),j),
\end{equation}
where $r\equiv p, j\equiv i\pmod{q}$ and $0\leq r, j<q$.  See
\cite{OSzRationalSurgeries} for details.

\subsection[L--space surgeries on knots]{$L$--space surgeries on knots} A rational homology
three-sphere $Y$ is called an {\em $L$--space\/} if $HF^+(Y)$ has no torsion and
the map $HF^\infty(Y)\rightarrow HF^+(Y)$ is surjective, or equivalently, if
$\smash{\widehat{HF}}(Y,\fraks)\cong\bbz$ for each \spinc structure $\fraks$ on $Y$.
The Heegaard Floer homology of an $L$--space is determined by its correction
terms. The following proposition is implied in \cite{OSzRationalSurgeries} (see
\cite{KMOSz} for results in the Seiberg--Witten gauge theory).

\begin{prop}\label{LSpaceSurgeryIsAHalfLine} Suppose $K$ is a nontrivial knot
in $S^3$ and $S^3_r(K)$ is an $L$--space with $r>0$. Then $S^3_s(K)$ is an
$L$--space for any rational number $s\geq r$.
\end{prop}

In particular, this implies that there will be no $L$--space surgeries on a
nontrivial amphichiral knot in $S^3$. This proposition is an ingredient to get
\fullref{hfkofcosmeticknot}.

\section{An equation in cyclotomic numbers}\label{cyclotomic}

We will need the following result regarding cyclotomic numbers in \fullref{cosmetic}. Note that the meaning of letters in this section are different
from other sections because of different conventions in number theory and
topology.

\begin{thm}\label{numbertheorem} Let $m$ be an integer at least $3$ and let $r$ and
$r^\prime$ be integers such that $(r,m)=(r^\prime,m)=1$. Then there exists an integer $k$ such that $(k,m)=1$ and
\begin{equation}\label{meqn1}
\babs{\Big}{\frac{(\xi-1+\xi^{-1})}{(\xi-1)(\xi^r-1)}}=
\babs{\Big}{\frac{(\xi^k-1+\xi^{-k})}{(\xi^k-1)(\xi^{r^\prime k}-1)}}
\end{equation}
for any $m$--th root of unity $\xi\neq1$ if and only if one of the following
holds.
\begin{enumerate}
\item $r\equiv\pm r^\prime\pmod{m}$.
\item $m=12$ and $\{\pm r, \pm r^\prime\}=\{\pm1, \pm5\}\pmod{12}$.
\item $m=18e+9$ and $\{\pm r,\pm r^\prime\}=\{\pm(6e+1),\pm(6e+5)\}\pmod{m}$
for some nonnegative integer $e$.
\end{enumerate}
\end{thm}

Let $\zeta=\zeta_m=e^{2\pi i/m}$. For $x=1,2,\cdots,m-1$, write
$A_x=\ln\abs{\zeta^x-1}$. Then when $\xi^6\neq1$, we can translate equation
\eqref{meqn1} into the following multiplicative relation of cyclotomic units:
\begin{equation}\label{meqn3}
\babs{\Big}{\frac{(\xi^6-1)}{(\xi^2-1)(\xi^3-1)(\xi^r-1)}}=
\babs{\Big}{\frac{(\xi^{6k}-1)}{(\xi^{2k}-1)(\xi^{3k}-1)(\xi^{r^\prime k}-1)}}
\end{equation}
In particular, when $\xi=\zeta$, in Ennola's notation, by taking
the logarithm this becomes
\begin{equation}\label{meqn2}
R:=A_6-A_2-A_3-A_r-\left(A_{6k}-A_{2k}-A_{3k}-A_{r^\prime k}\right)=0.
\end{equation}
The ``if" part of the proof is easy and left to the reader. The machinery in
the proof of the other direction are from Franz \cite{Franz} and Dvornicich
\cite{Dvornicich}. Also, when $m\leq8$, we may compute the torsion directly to
prove the theorem. So in what follows, we assume $m>8$. Since for
$m\equiv2\pmod{4}$ there are no primitive characters, we will deal separately
with the case $m\not\equiv2\pmod{4}$ in \fullref{primitive} and with
the case $m\equiv2\pmod{4}$ in \fullref{noprimitive}.

\subsection{Preliminaries in number theory}
We refer to Washington \cite{Washington} or Apostol \cite{Apostol} (which is more fundamental) for
basic notation.

In 1935, Franz used the following lemma to classify lens spaces
(in any dimension).

\begin{thm}[Franz Independence Lemma \cite{Franz}] Fix a natural number $m$ and
let $S=\{j\in\bbz/m\bbz\ |\ (j, m)=1 \}$. Suppose that $\{a_k\ |\ k\in S\}$ is
a set of integers satisfying the following.
\begin{enumerate}
\renewcommand{\labelenumi}{(\alph{enumi})}
\item $\sum_{k\in S}a_k=0$.
\item $a_j=a_{-j}$ for all $j\in S$.
\item $\prod_{j\in S}(\xi^j-1)^{a_j}=1$ for every $m$--th root of unity $\xi\neq1$.
\end{enumerate}
Then $a_j=0$ for all $j\in S$.
\end{thm}

This lemma says that there are no multiplicative relations for the set
$\{\zeta_m^j\ |\ j\in S\}$ where $\zeta_m=e^{2\pi i/m}$. For a proof for
topologists, we refer to the lecture notes of de Rham
\cite[Chapter 1]{deRham}.

For an integer $m\geq2$. Let $\zeta=\zeta_m=e^{2\pi i/m}$ and $A_x=\ln\abs{\zeta^x-1}$ for $x\in\bbz/m\bbz\setminus\{0\}$. Let $R=\sum_{x=1}^{m-1}C_xA_x$ be a linear combination of $A_x$ with integer coefficients. A mod ${n}$ character is a multiplicative homomorphism $\chi\co (\bbz/n\bbz)^*\rightarrow \bbc^*$ and we let $\chi(i)=0$ if $(i,n)>1$. If $n|m$, then a mod $n$ character induces a  mod $m$ character by composing the homomorphism $(\bbz/n\bbz)^*\rightarrow(\bbz/m\bbz)^*$. The conductor of a mod $m$ character $\chi$ is the smallest $n$ such that $\chi$ is induced by a mod $n$ character. For a character $\chi \pmod{m}$ with conductor $f>1$ and integer
$d\geq1$ such that $f\mid d$ and $d\mid m$, we denote
$$T(\chi, d, R)=\underset{(x,d)=1}{\sum_{x=1}^{d-1}}\chi(x)C_{(m/d)x}$$
$$Y(\chi, R)=\underset{f\mid d,\, d\mid m}{\sum_d}\frac{1}{\phi(d)}\prod_{p\mid d}
(1-\wbar{\chi}(p))T(\chi, d, R).\leqno{\hbox{and}}$$ For $p\mid m$, define $\gamma_p$ by
$p^{\gamma_p}\mid m$ and $p^{\gamma_p+1}\nmid m$ and denote
$$Y_p(R)=\sum_{x=1}^{p^{\gamma_p}-1}(x, p^{\gamma_p}) C_{(m/p^{\gamma_p})x}.$$
Ennola gave the following characterization of relations in
cyclotomic numbers:

\begin{thm}[Ennola \cite{Ennola}] With all notation as above, $R=0$ if and only if
the following two conditions hold:
\begin{align}\label{ennolarelation1}
Y(\chi, R)&=0\quad\text{for every even character $\chi\neq\chi_1$}\\
Y_p(R)&=0\quad\text{for every prime number $p$ dividing $m$}\notag
\end{align}
where $\chi_1$ is the principle character defined by $\chi_1(x)=1$ for
$(x,m)=1$ (and $\chi_1(x)=0$ if $(x,m)>1$).
\end{thm}

In our application of this theorem, the second relations that
$Y_p(R)=0$ for every prime number $p$ dividing $m$
are
automatically satisfied. We shall prove our results as consequence of the
relations \eqref{ennolarelation1}.

Dvornicich gave the following propositions about
characters (in order to solve his equation in cyclotomic numbers):

\begin{prop}{\rm (Dvornicich \cite[Corollary 1]{Dvornicich})}\qua\label{dvo1} Let $m\geq2$. If $k$
satisfies $\chi(k)=1$ for all even primitive characters of a maximal conductor
$C$, where $C=m/2$ if $m\equiv 2\pmod{4}$ and $C=m$ otherwise, then
\begin{enumerate}
\renewcommand{\labelenumi}{(\alph{enumi})}
\item $k\equiv \pm1,\pm i\pmod{m}$ if $m=2^a\cdot3, a\geq3$, where
$i=1-2^{a-1}$ if $a$ is even and $i=1+2^{a-1}$ when $a$ is odd.
\item $k\equiv\pm1\pmod{m/2}$ if $m\equiv2\pmod{4}$.
\item $k\equiv\pm1\pmod{m}$ otherwise.
\end{enumerate}
\end{prop}

\begin{prop}{\rm (Dvornicich \cite[Lemma 3]{Dvornicich})}\qua\label{dvo2} If $m\geq3$ is odd. Let
$G$ be the group of all even Dirichlet characters mod $m$. Then the set of
Dirichlet even characters of maximal conductor is not contained in the union of
two maximal subgroups of $G$, unless $m=3m_1m_2$ with $m_1,m_2>1$ and
$(3,m_1)=(3,m_2)=(m_1,m_2)=1$.
\end{prop}

The following lemma is told to me by Bjorn Poonen in private communication.
This lemma says that if four roots of unity add up to zero, then they must
cancel in pairs.

\begin{lemma}\label{fourroots} Suppose $\xi_i, 1\leq i\leq4$, are all roots of
unity and $\xi_1+\xi_2=\xi_3+\xi_4$.  Then $\xi_1+\xi_2=0$, $\xi_1=\xi_3$ or
$\xi_1=\xi_4$.
\end{lemma}

\begin{proof} Suppose that $\xi_1+\xi_2\neq0$. Considering the angle and
norm of $\xi_1+\xi_2$ and $\xi_3+\xi_4$ will give the lemma.\end{proof}

\subsection[Proof of theorem when m is not congruent to 2 mod 4]{Proof of theorem when $\bf m\not\equiv 2\text{ (mod 4)}$}
\label{primitive} Pick an integer $k^\prime$ such that $r^\prime\equiv
rk^\prime\pmod{m}$.

{\bf Case 1}\qua $(m, 6)=1$.

In this case, we have $\xi^6\neq1$ for all $m$--th roots of unity $\xi\neq1$.
Hence equation \eqref{meqn3} holds for all $m$--th roots of unity $\xi\neq1$ and
is nonzero. So the Franz Independence Lemma applies and we get
$$\{\pm6, \pm2k,\pm3k,\pm r^\prime k\}=\{\pm6k, \pm2,\pm3,\pm r\}\pmod{m}.$$
By congruence arguments, we can show that $r\equiv\pm r^\prime\pmod{m}$.

{\bf Case 2}\qua $(m, 12)=4$.

Let $b$ be the mod $m$ solution of the congruence equation
$r\equiv3bk\pmod{m}$. Consider all relations \eqref{ennolarelation1} relative
to characters of maximal conductor. They have the form
$$\chi(3)+\chi(r)-\chi(3k)-\chi(rk^\prime k)=0.$$
(All other terms are zero since they are even and not coprime to $m$.)
By \fullref{fourroots} we have the following three possibilities.
\begin{enumerate}
\item $\chi(3)+\chi(r)=\chi(3k)+\chi(r k^\prime k)=0$. Since $\chi(k)\neq0$, it
follows that $\chi(r)=\chi(rk^\prime)$ and hence $\chi(k^\prime)=1$.
\item $\chi(3)-\chi(3k)=\chi(r)-\chi(rk^\prime k)=0$. We get $\chi(k)=1$ and
hence $\chi(k^\prime)=1$.
\item $\chi(3)-\chi(rk^\prime k)=\chi(r)-\chi(3k)=0$. So $\chi(3bk)-\chi(3k)=0$
and $\chi(b)=1$.
\end{enumerate}
In all cases, we get
$$(\chi(k^\prime)-1)(\chi(b)-1)=0.$$
Suppose $k^\prime\not\equiv\pm1,b\not\equiv\pm1\pmod{m}$. Let $G_1$, $G_2$ be
the subgroup of $G$ defined by $\chi(k^\prime)=1$ and $\chi(b)=1$. By
\fullref{dvo1}, $G_1, G_2$ are proper subgroups of $G$. Now by
\fullref{dvo2}, we have the equalities $m=3m_1m_2$ with $m_1,m_2>1$ and
$(3,m_1)=(3,m_2)=(m_1,m_2)=1$, which is impossible since $(m,3)=1$. Therefore
either $k^\prime\equiv\pm1\pmod{m}$ or $b\equiv\pm1\pmod{m}$.

If $k^\prime\equiv\pm1\pmod{m}$, we then have $r\equiv\pm r^\prime\pmod{m}$.

If $b\equiv\pm1\pmod{m}$, then $r\equiv\pm 3k\pmod{m}$ and our equation
\eqref{meqn2} is reduced to
\begin{equation}\label{case2-1}
A_6-A_2-A_3-A_{6k}+A_{2k}+A_{r^\prime k}=0.
\end{equation}
Now for any even character of maximal conductor,
\eqref{ennolarelation1} becomes
$$-\chi(3)+\chi(r^\prime k)=0.$$
Again by \fullref{dvo1}, we get $r^\prime k\equiv\pm 3\pmod{m}$. This
reduces our equation further to
$$A_6-A_2-A_{6k}+A_{2k}=0.$$
Rewrite this equation to get
$$\abs{\zeta^2+1+\zeta^{-2}}=\abs{\zeta^{2k}+1+\zeta^{-2k}} \ \ \text{and}\ \
\abs{1+2\cos({4\pi/m})}=\abs{1+2\cos(4k\pi/m)}.$$ Note that since $m>8$, we have
$\cos(4\pi/m)>0$; hence $\cos(4k\pi/m)>0$ and we have
$$\zeta^2+1+\zeta^{-2}=\zeta^{2k}+1+\zeta^{-2k}.$$
By \fullref{fourroots}, $2k\equiv\pm2\pmod{m}$. So either
$k\equiv\pm1\pmod{m}$ or $k\equiv m_1\pm1\pmod{m}$, where $m_1=m/2$. If
$k\equiv\pm1\pmod{m}$, we get $r=\pm3, r^\prime=\pm3\pmod{m}$. If $k\equiv
m_1\pm1\pmod{m}$, then $r\equiv\pm(m_1+3)$ and $r^\prime\equiv
\pm(m_1+3)\pmod{m}$. Therefore $r\equiv\pm r^\prime\pmod{m}$.

{\bf Case 3}\qua $(m, 6)=3$.

Let $b$ be the mod $m$ solution of the congruence equation
$r\equiv2bk\pmod{m}$. Consider all relations \eqref{ennolarelation1} relative
to characters of maximal conductor. They have the form
$$\chi(2)+\chi(r)-\chi(2k)-\chi(rk^\prime k)=0.$$
\newpage
By \fullref{fourroots}, we have
the following three possibilities.
\begin{enumerate}
\item $\chi(2)+\chi(r)=\chi(2k)+\chi(r k^\prime k)=0$.
Hence $\chi(k)\left(\chi(2)+\chi(rk^\prime)\right)=0$ and it follows that
$\chi(r)=\chi(rk^\prime)$ and therefore $\chi(k^\prime)=1$.
\item $\chi(2)-\chi(2k)=\chi(r)-\chi(rk^\prime k)=0$. We get $\chi(k)=1$ and
hence $\chi(k^\prime)=1$.
\item $\chi(2)-\chi(rk^\prime k)=\chi(r)-\chi(2k)=0$. So $\chi(2bk)-\chi(2k)=0$
and $\chi(b)=1$.
\end{enumerate}
In all three cases, we get
\begin{equation}\label{case3}
(\chi(k^\prime)-1)(\chi(b)-1)=0.
\end{equation}
Suppose $k^\prime\not\equiv\pm1,b\not\equiv\pm1\pmod{m}$. Let $G_1$, $G_2$ be
the subgroup of $G$ defined by $\chi(k^\prime)=1$ and $\chi(b)=1$. By
\fullref{dvo1}, $G_1, G_2$ are proper subgroups of $G$.
Thus we have $m=3m_1m_2$ with $m_1,m_2>1$ and
$(3,m_1)=(3,m_2)=(m_1,m_2)=1$ by
\fullref{dvo2}.

So there are three possibilities.

(1)\qua If $k^\prime\equiv\pm1\pmod{m}$, we get $r\equiv\pm r^\prime\pmod{m}$.

(2)\qua If $b\equiv\pm1\pmod{m}$, ie, $r\equiv\pm 2k\pmod{m}$, equation
\eqref{meqn2} is reduced to
\begin{equation*}
A_6-A_2-A_3-A_{6k}+A_{3k}+A_{r^\prime k}=0.
\end{equation*}
Now for any even character of maximal conductor,
\eqref{ennolarelation1} becomes
$$-\chi(2)+\chi(r^\prime k)=0.$$
Again by \fullref{dvo1}, we get $r^\prime k\equiv\pm 2\pmod{m}$. This
reduces our equation further to
$$A_6-A_3-A_{6k}+A_{3k}=0.$$
This gives $\abs{\zeta^3+1}=\abs{\zeta^{3k}+1}$. Hence $3k\equiv\pm 3\pmod{m}$
and we get the congruence system
\begin{equation}\label{case3-2}
3k\equiv\pm3,\quad r\equiv\pm2k,\quad r^\prime k\equiv \pm 2\quad \pmod{m}.
\end{equation}
Using congruence arguments, we can show that either $r\equiv\pm
r^\prime\pmod{m}$, or,\ for some nonnegative integer $e$,
$$m=18 e+9\quad\text{and}\quad\{\pm r, \pm r^\prime\}=\{\pm(6e+1), \pm(6e+5)\}.$$
(3)\qua If $m=3m_1m_2$ with $m_1,m_2>1$ and $(3,m_1)=(3,m_2)=(m_1,m_2)=1$. In
equation \eqref{meqn3}, let $\xi=\zeta_m^{3a}$ with $m_1m_2\nmid a$, and we
find that equation \eqref{meqn3} holds for all $m_1m_2$--th root of unity
$\xi\neq1$ and $(m_1m_2,6)=1$. So by case 1, we get $r\equiv\pm
r^\prime\pmod{m_1m_2}$ and $k\equiv\pm1\pmod{m_1m_2}$. This implies
$3k\equiv\pm3, 6k\equiv\pm6\pmod{m}$.

So equation \eqref{meqn3} reduces to
$$\abs{(\xi^2-1)(\xi^r-1)}=\abs{(\xi^{2k}-1)(\xi^{r^\prime k}-1)}.$$
Now the Franz Independence Lemma applies and we get
$$\{\pm2,\pm r\}=\{\pm2k, \pm r^\prime k\}\quad \text{and}\quad3k\equiv\pm3\pmod{m}.$$
If $2k\equiv\pm2\pmod{m}$, then we have $k\equiv \pm1, \pm 5\pmod{m}$. If
$k\equiv\pm5\pmod{m}$, we get $10\equiv\pm 2, 15\equiv\pm3\pmod{m}$, which
implies $m=12$ and which is impossible since $m$ is odd. So we have
$k\equiv\pm1\pmod{m}$ and $r\equiv\pm r^\prime\pmod{m}$.

If $2k\equiv\pm r\pmod{m}$, then we get the same congruence system
\eqref{case3-2}. So $r\equiv\pm r^\prime\pmod{m}$ since $m$ can not be of the
form $18e+9$.

This ends the proof of the theorem in Case 3.

{\bf Case 4}\qua $(m, 12)=12$.

For all characters of maximal conductor, Ennola's relations
\eqref{ennolarelation1} have the form
$$\chi(r)-\chi(r^\prime k)=0, \quad \text{ie,}\quad \chi(r)-\chi(rk^\prime k)=0.$$
So we get $\chi(k^\prime k)=1$. Now by \fullref{dvo1}, we have one of
the following.
\begin{enumerate}
\item $k^\prime k\equiv\pm 1$ if $m$ is not of the form $2^a\cdot3$, $a\geq3$.
\item $k^\prime k\equiv\pm 1,\pm i\pmod{m}$ if $m=2^a\cdot3$, $a\geq3$. Here
$i=1-2^{a-1}$ if $2\mid a$ and $i=1+2^{a-1}$ if $2\nmid a$.
\end{enumerate}

{\bf Subcase 4.1}\qua $m\neq 2^a\cdot3$ for any $a\geq 3$.

In this subcase, we will get $r\equiv \pm r^\prime k\pmod{m}$. So our equation
\eqref{meqn2} reduces to
\begin{equation}\label{case41}
R^\prime:=A_6-A_2-A_3-\left(A_{6k}-A_{2k}-A_{3k}\right)=0.
\end{equation}
Write $m=12m_1$. Consider even characters $\chi$ of conductor $f_1=4m_1$.

Let $d$ be an integer such that $f_1|d|m$. Since $m/f_1=3$, $m/d=1,3$. Let
$d_1=12m_1$ and $d_2=4m_1$. For $d_1=12m_1$, we have
$$T(\chi, d_1, R^\prime)=\underset{(x,m)=1}{\sum_{x=1}^{m-1}}\chi(x)C_x=0$$
since $(x,m)>1$ for every $C_x\neq0$. For $d_2$, we get
$$T(\chi, d_2, R^\prime)=-\chi(1)+\chi(k)=\chi(k)-1.$$
\newpage
So \eqref{ennolarelation1} for such characters becomes
$$\frac1{\varphi(f_1)}(\chi(k)-1)=0.$$
So $\chi(k)=1$ and by \fullref{dvo1} we get $k\equiv\pm1\mod{f_1}$, or
$k\equiv\pm i\pmod{f_1}$ if $f_1=2^a\cdot3, a\geq3$ (equivalently
$m=2^a\cdot9$).

If $k\equiv\pm i\pmod{f_1}$, then $m=2^a\cdot9$, $a\geq3$, and we have
\begin{equation}\label{case41-2}
\begin{split}
k\equiv\pm(1-2^{a-1})\pmod{2^a}&\quad\text{ if $2\mid a$ }\\
k\equiv\pm(1+2^{a-1})\pmod{2^a}&\quad\text{ if $2\nmid a$.}
\end{split}
\end{equation}
For any even character $\chi$ of conductor $f=2^a$,
we have $m/f=9$ and $d$ could be $d_1=m$, $d_2=m/3$ and $d_3=m/9$. For our
equation \eqref{case41}, we have
$$T(\chi, d_1, R^\prime)=\underset{(x,m)=0}{\sum_{x=1}^{m-1}}\chi(x)C_x=0,\quad
T(\chi, d_3, R^\prime)=\underset{(x, 2^a)=0}{\sum_{x=1}^{d_3-1}}
\chi(x)C_{9x}=0,$$
$$T(\chi, d_2, R^\prime)=\underset{(x,
2^a\cdot3)=0}{\sum_{x=1}^{d_2-1}} \chi(x)C_{3x}=-\chi(1)+\chi(k).$$ So $Y(\chi,
R^\prime)=0$ gives
$-\chi(1)+\chi(k)=0$.
By \fullref{dvo1}, we get $k\equiv\pm1\pmod{2^a}$. This is a
contradiction to \eqref{case41-2} since $a\geq3$. So we must have
$k\equiv\pm1\pmod{f_1}$.

If $k\equiv\pm1\pmod{f_1}$, then we have $6k\equiv\pm6,3k\equiv\pm3\pmod{m}$.
And equation \eqref{case41} implies $A_2=A_{2k}$. So we have the system
$$\begin{cases}2k\equiv\pm2\pmod{m},\\3k\equiv\pm3\pmod{m}.\end{cases}$$
So $k\equiv\pm1,\pm5\pmod{m}$. If $k\equiv\pm1\pmod{m}$, then
$k^\prime\equiv\pm1\pmod{m}$ and hence $r\equiv\pm r^\prime\pmod{m}$. If
$k\equiv\pm5\pmod{m}$, then $10\equiv\pm2\pmod{m}$ and hence $m=12$ since
$12|m$. By an easy computation, we can conclude $\{\pm r,\pm
r^\prime\}=\{\pm1,\pm5\}$.

This ends the proof of theorem in subcase 4.1.

{\bf Subcase 4.2}\qua $m=2^a\cdot3$, $a\geq 3$.

Let us first suppose $k^\prime k\equiv\pm i\pmod{m}$. If $a$ is even, then we
get $k^\prime k\equiv \pm(1-2^{a-1})\pmod{m}$. If $a$ is odd, then $k^\prime
k\equiv\pm(1+ 2^{a-1})\pmod{m}$.

Let $\chi_0$ be a generator of the even characters mod ${2^a}$; we
have $\chi_0(\pm k^\prime k)=-1$ and the relation $Y(\chi_0, R)=0$
becomes
$$\frac{1}{2^{a-1}}\underset{(x,2^a)=1}{\sum_{x=1}^{2^a-1}}\chi_0(x)C_{3x} +
\frac{1}{2\cdot 2^{a-1}}(1-\wbar{\chi}_0(3))\underset{(x,2^a\cdot3)=1}
{\sum_{x=1}^{2^a\cdot3-1}}\chi_0(x)C_{x}=0.$$ And for our relation $R$, we read
$$2(-1+\chi_0(k))+(1-\wbar{\chi}_0(3))(-\chi_0(r)+\chi_0(rk^\prime k))=0.$$
Since $\chi_0(k^\prime k)=-1$, we get
$$\chi_0(k)-1=\chi_0(r)(1-\wbar{\chi}_0(3)).$$
This gives a $4$--term relation among roots of unity. By \fullref{fourroots}, since $\chi_0(3)\neq1$, we have the following two
possibilities.
\begin{enumerate}
\item $\chi_0(k)=\chi_0(r)$, $1=\chi_0(r)\wbar{\chi}_0(3)$. Hence we have
$\chi_0(k)=\chi_0(r)=\chi_0(3)$.
\item $\chi_0(k)=-\chi_0(r)\wbar{\chi}_0(3)$, $\chi_0(r)=-1$. Let
$\bar{k}, \bar{k}^{\prime}$ be integers where
$\bar{k}^{\prime}k^\prime\equiv1$ and $\bar{k}k\equiv1\pmod{m}$. Using
$\chi_0(k^\prime k)=-1$, we may conclude that these two equations will be
equivalent to $\chi_0(\bar{k})=\chi_0(\bar{k}^{\prime})=\chi_0(3)$.
\end{enumerate}

By symmetry, we may just consider one possibility. Suppose
$\chi_0(k)=\chi_0(r)=\chi_0(3)$.

Consider the character $\chi_1=\chi_0^2$. It is a character of conductor
$2^{a-1}=m/6$. There are only four possibilities: $d_1=m$, $d_2=m/2=3\cdot
2^{a-1}$, $d_3=m/3=2^a$ and $d_4=m/6=2^{a-1}$. Since $\chi_0(k^\prime
k)=-1$ and $\chi_0(k)=\chi(3)$, we have
\begin{align*}T(\chi_1, d_1, R)&=-\chi_1(r)+\chi_1(rk^\prime k)=\chi_1(r)\left(
\chi_0^2(k^\prime k)-1\right)=0,\\
T(\chi_1, d_2, R)&=\underset{(x,3\cdot 2^{a-1})=1}{\sum_{x=1}^{d_2-1}}\chi_1(x)C_{2x}
=-\chi_1(1)+\chi_1(k)=\chi_0^2(3)-1,\\
T(\chi_1, d_3, R)&=\chi_0^2(3)-1 \quad\text{and}\quad T(\chi_1, d_4, R)=1-\chi_0^2(3).\end{align*}
Therefore, since $\chi_1(2)=0$, the relation \eqref{ennolarelation1} for the character $\chi_1$
becomes
\begin{align*}Y(\chi_1, R)=&\frac{1}{2\cdot
2^{a-2}}(1-\wbar{\chi}_0^2(3))(\chi_0^2(3)-1)
+\frac{1}{2^{a-1}}(\chi_0^2(3)-1) + \frac{1}{2^{a-2}}(1-\chi_0^2(3))\\
=&\frac{(\chi_0^2(3)-1)}{2^{a-1}}\left((1-\wbar{\chi}_0^2(3))+1-2\right)
=\frac{\wbar{\chi}_0^2(3)(\chi_0^2(3)-1)}{2^{a-1}}.
\end{align*}
\newpage
For $a=3$, we may prove the theorem by direct computation. And when $\chi_0^2(3)\neq 1$, this is a contradiction. Hence $k^\prime k\equiv
\pm1\pmod{m}$.

When $k^\prime k\equiv\pm1\pmod{m}$, we reduce our equation to
$$R^\prime:=A_6-A_2-A_3-\left(A_{6k}-A_{2k}-A_{3k}\right)=0.$$
For the character $\chi_0$, we get
$$Y(\chi_0, R^\prime)=\frac1{2^{a-1}}(-1+\chi_0(k))=0,$$
and hence $\chi_0(k)=1$. Therefore $k\equiv\pm1\pmod{2^a}$. This implies
$3k\equiv\pm3, 6k\equiv\pm6\pmod{m}$ and $R^\prime=0$ will reduce to
$A_2=A_{2k}$. So we get $2k\equiv \pm2,\pm( m/2+2)\pmod{m}$. Combining
$k\equiv\pm1\pmod{2^a}$, we get $k\equiv\pm1\pmod{m}$. Hence
$k^\prime\equiv\pm1 \pmod{m}$, which implies $r\equiv\pm r^\prime\pmod{m}$.
This ends the proof of the theorem in subcase 4.2.

\subsection[Proof of theorem when m is congruent to 2 mod 4]{Proof of theorem when $\bf m\equiv 2\text{ (mod 4)}$}
\label{noprimitive} There is no primitive character with modulus $m$. Write
$m=2m_1$, where $m_1=2a+1$ is odd. For any $m_1$--th root of unity $\xi$, $-\xi$
is a $m$--th root of unity, hence we have
\begin{equation}\label{meqn2-*-1}\babs{\Big}{\frac{\xi+1+\xi^{-1}}{(\xi+1)(\xi^r+1)}}=
\babs{\Big}{\frac{\xi^k+1+\xi^{-k}}{(\xi^k+1)(\xi^{r^\prime k}+1)}}.
\end{equation}
When
$\xi\neq1$, this could be transformed into
\begin{equation}\label{meqn2-*-2}\babs{\Big}{\frac{(\xi^3-1)(\xi^{r}-1)}{(\xi^2-1)(\xi^{2r}-1)}}=
\babs{\Big}{\frac{(\xi^{3k}-1)(\xi^{r^\prime k}-1)}{(\xi^{2k}-1)(\xi^{2r^\prime
k}-1)}}.
\end{equation}
In particular, this holds for $\xi=\zeta_{m_1}$. So in Ennola's notation, we have
the mod $m_1$ relation
$$R_1:=A_{3}+A_{r}-A_{2}-A_{2r}-(A_{3k}+A_{r^\prime k}-A_{2k}-A_{2r^\prime k})=0.$$
Note that $r, k, r^\prime$ are all coprime to $2$ and $m_1$.

{\bf Case 1}\qua $(m_1, 3)=1$.

In this case, equation \eqref{meqn2-*-2} is
$$\babs{\Big}{\frac{(\xi^3-1)(\xi^{r}-1)}{(\xi^2-1)(\xi^{2r}-1)}}=
\babs{\Big}{\frac{(\xi^{3k}-1)(\xi^{r^\prime k}-1)}{(\xi^{2k}-1)(\xi^{2r^\prime
k}-1)}}\neq0,$$
\newpage
and holds for any $m_1$--th root of unity $\xi\neq1$. So the Franz
Independence Lemma implies that
$$\{\pm 3,\pm r,\pm2k,\pm2r^\prime k\}=\{\pm 3k,\pm r^\prime k, \pm2,\pm 2r\}\pmod{m_1}.$$
By a congruence argument, we get $r\equiv\pm r^\prime\pmod{m_1}$ or
$$r\equiv\pm2, k\equiv \pm2, r^\prime\equiv\pm \frac{m_1+3}{2}\pmod{m_1}$$
or a symmetric case. In the first case, we get $r\equiv
\pm r^\prime\pmod{m}$ since $r\equiv\pm r^\prime\pmod{2}$.
In the second case, we get $r\equiv \pm k$ and $r^\prime k\equiv \pm3\pmod{m_1}$. Note that equation \eqref{meqn1} also holds for any $m_1$--th root $\xi$, so we have
$$\babs{\Big}{\frac{\xi-1+\xi^{-1}}{\xi-1}}=\babs{\Big}{\frac{\xi^2-1+\xi^{-2}}{\xi^3-1}}.$$
When $\xi\neq-1$, this could be simplified as
$\abs{\xi^2+1+\xi^{-2}}=\abs{\xi^2-1+\xi^{-2}}$, or equivalently, $\abs{2{\rm Re}(\xi^2)+1}= \abs{2{\rm Re}(\xi^2)-1}$. We hence get ${\rm Re}(\xi^2)=0$ for any $m_1$--th root of unity other than $\pm1$, which is a contradiction. So we must have $r\equiv\pm r^\prime\pmod{m}$.

{\bf Case 2}\qua $(m_1, 3)=3$.

For any even character $\chi$ of maximal conductor,
\eqref{ennolarelation1} has the form
\begin{equation}\label{2-2-1}
\chi(r)-\chi(2)-\chi(2r)-\chi(r^\prime k)+\chi(2k)+\chi(2 r^\prime k)=0.
\end{equation}
In equation \eqref{meqn3}, let $\xi=\zeta_m^2=\zeta_{m_1}$. Then no term in
equation \eqref{meqn3} will be zero (note that we suppose $m>8$). And we will
get equation \eqref{meqn2} as a mod $m_1$ relation. So for any even character
$\chi$ of maximal conductor, we get
$$-\chi(2)-\chi(r)+\chi(2k)+\chi(r^\prime k)=0.$$
Substitute $\chi(r^\prime k)$ in \eqref{2-2-1} by $\chi(2)+\chi(r)-\chi(2k)$ to
get
$$\chi(2)(2-\chi(2))(\chi(k)-1)=0.$$
Hence $\chi(k)=1$ and by \fullref{dvo1}, we get
$k\equiv\pm1\pmod{m_1}$. Since $(m,k)=1$, we conclude that
$k\equiv\pm1\pmod{m}$. Now equation \eqref{meqn2-*-1} reduces to
$$\abs{\xi^r+1}=\abs{\xi^{r^\prime}+1},$$
and we get $r\equiv\pm r^\prime\pmod{m}$.

This ends the proof of our theorem.

\section{Cosmetic surgeries on genus one knots}\label{cosmetic}

In this section, we prove \fullref{genusonecosmetic}.

Our focus will be on knots having the same knot Floer homology as the right-handed trefoil. For a knot $K$ in $S^3$, we call any knot having the same knot
Floer homology as $K$ a {\em fake $K$ knot\/}. In particular, $K$ is itself a
fake $K$ knot.

Using the notation of \fullref{hfhprelim}, for any fake right-handed trefoil
knot $K$, we have $H_*(A_s^+)\cong H_*(B^+)\cong\calt^+$, where $\calt^+$ is
the $\bbz[U]$--module $\bbz[U,U^{-1}]/(U\cdot\bbz[U])$. The two maps on
homology (still denoted by $h^+$ and $v^+$) induced by $h^+\co A^+_s\rightarrow
B^+$ and $v^+\co A^+_s\rightarrow B^+$ are identified with multiplication by
powers of $U$ as follows.
\begin{equation}\label{trefoilmaps}
h^+_s=\begin{cases}1\cdot&\text{if $s<0$}\\
U\cdot&\text{if $s=0$}\\ U^s\cdot&\text{if $s>0$}\end{cases}\qquad
v^+_s=\begin{cases}1\cdot&\text{if $s>0$}\\U\cdot&\text{if $s=0$}\\
U^{-s}\cdot&\text{if $s<0$}\end{cases}
\end{equation}
Observe that $h^+_s=v^+_{-s}$, and we get a symmetry of the chain complex
$\bbx^+_{\smash{p/q}}$ given by
\begin{equation}\label{trefoilsymmetry}
(t,x)\longleftrightarrow(q-t-1,x),\quad x\in\calt^+.
\end{equation}
As a motivational application of the rational surgery (to our propositions), we
give the chain complex $\bbx^+_{5/3}(K)$ illustrated as follows.
$$
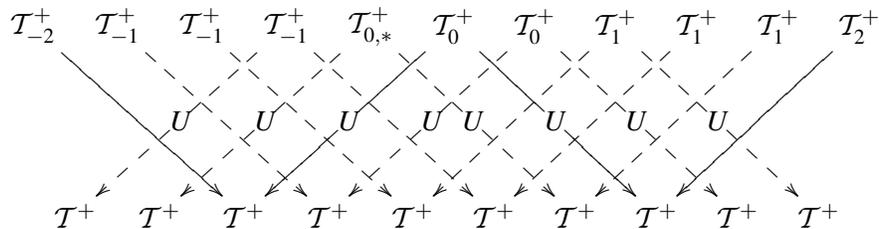
\xymatrix@C=-6.5pt{ \calt^+_{-2}  \ar[ddrrrrr] && \calt^+_{-1}  \ar@{-->}[ddrrrrr]
&& \calt^+_{-1} \ar@{-->}[ddrrrrr] && \calt^+_{-1}\ar@{-->}[ddlllll]|{\dps{U}}
\ar@{-->}[ddrrrrr] && \calt^+_{0,*} \ar@{-->}[ddlllll]|{\dps{U}}
\ar@{-->}[ddrrrrr]|{\dps{U}} && \calt^+_0 \ar[ddlllll]|{\dps{U}}
\ar[ddrrrrr]|{\dps{U}} &&\calt^+_0 \ar@{-->}[ddlllll]|{\dps{U}}
\ar@{-->}[ddrrrrr]|{\dps{U}} && \calt^+_1 \ar@{-->}[ddlllll]
\ar@{-->}[ddrrrrr]|{\dps{U}} && \calt^+_1 \ar@{-->}[ddlllll]&& \calt^+_1
\ar@{-->}[ddlllll] &&\calt^+_2 \ar[ddlllll]\\
&&&&&&&&&&&&\\&\calt^+ && \calt^+ && \calt^+ && \calt^+ && \calt^+ &&\calt^+ &&
\calt^+ && \calt^+ && \calt^+ && \calt^+}$$ In this diagram, the one with a
star subscript corresponds to $(0, A^+_0)$, while the number in the subscript
means that $\calt^+_i$ corresponds to $A^+_i$. Also the right arrows denote
$h^+$ while the left arrows denote $v^+$. The label at the middle of an arrow
gives the map in the sense of \eqref{trefoilmaps}, and those unlabeled ones are
identities.
\newpage
Similarly, the chain complex $\bbx^+_{5/4}(K)$ could be illustrated as follows.
$$
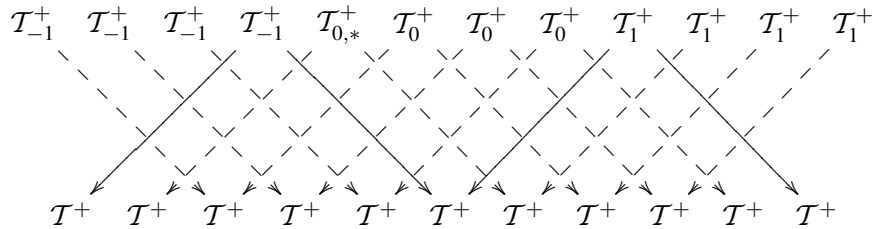
\xymatrix@C=-8pt{  \calt^+_{-1}  \ar@{-->}[ddrrrrr] &&  \calt^+_{-1}
\ar@{-->}[ddrrrrr] && \calt^+_{-1} \ar@{-->}[ddrrrrr] && \calt^+_{-1}
\ar@{->}[ddlllll] \ar@{->}[ddrrrrr] && \calt^+_{0,*}\ar@{-->}[ddlllll]
\ar@{-->}[ddrrrrr] && \calt^+_{0} \ar@{-->}[ddlllll] \ar@{-->}[ddrrrrr] &&
\calt^+_0 \ar@{-->}[ddlllll] \ar@{-->}[ddrrrrr] &&\calt^+_0 \ar@{-->}[ddlllll]
\ar@{-->}[ddrrrrr] && \calt^+_1 \ar@{->}[ddlllll] \ar@{->}[ddrrrrr] &&
\calt^+_1 \ar@{-->}[ddlllll]&& \calt^+_1
\ar@{-->}[ddlllll] &&\calt^+_1 \ar@{-->}[ddlllll]\\
&&&&&&&&&&&&\\&\calt^+ &&\calt^+ && \calt^+ && \calt^+ && \calt^+ && \calt^+
&&\calt^+ && \calt^+ && \calt^+ && \calt^+ && \calt^+}$$ Here we omit the
illustration of homomorphisms, but one could identify them using formula
\eqref{trefoilmaps}.

We now establish some propositions on the Heegaard Floer homology of a surgered
manifold on a fake right-handed trefoil knot.

\begin{prop}\label{NoNegativeOnRHT} Suppose $K$ is a fake right-handed trefoil
knot, then no negative surgery on $K$ yields an $L$--space. And symmetrically,
there are no positive surgeries on fake left-handed trefoil knots yielding
$L$--spaces.
\end{prop}
\begin{proof} This is an application of the rational surgery formula (\fullref{rationalsurgeryformula}).
\end{proof}

\begin{prop}\label{NoSmallLSpaceSurgery} For any nontrivial knot $K$ in $S^3$,
if $S^3_r(K)$ is an $L$--space, then $\abs{r}\geq1$.
\end{prop}

\begin{proof} This is a direct consequence of the rank formula of
$\widehat{HF}(S^3_r(K))$ \cite[Proposition 9.5]{OSzRationalSurgeries}.
\end{proof}

\begin{prop}\label{UnknotTrefoil} Let $p$ and $q$ be a pair of positive coprime
integers. Suppose that $p$ is odd and $q\leq p$. Let $U$ be the unknot and $K$
be a fake right-handed trefoil knot, then, for the spin \spinc structure
$\fraks_0$, we have
$$d(S^3_{\smash{p/q}}(K), \fraks_0)=\begin{cases}
d(S^3_{\smash{p/q}}(U),\fraks_0),&\text{ if $q$ is even},\\
d(S^3_{\smash{p/q}}(U),\fraks_0)-2,&\text{ if $q$ is odd}.\end{cases}$$
\end{prop}

\begin{proof} We apply the rational surgery formula for Heegaard Floer homology
(\fullref{rationalsurgeryformula}) here.

In view of the property of the correction terms
\eqref{CorrectionTermProperties}, the symmetry mentioned above implies that
the central subcomplex which is preserved under the symmetry
\eqref{trefoilsymmetry} has the same homology as that corresponding to the spin
\spinc structure $\fraks_0$. (In the diagrams illustrating the chain complexes
for $S^3_{5/3}(T)$ and $S^3_{5/4}(T)$, the two central complexes are
represented by solid lines.)

When $q$ is odd (see the complex for $S^3_{5/3}(K)$ as an example), this
central subcomplex has the following form:
$$
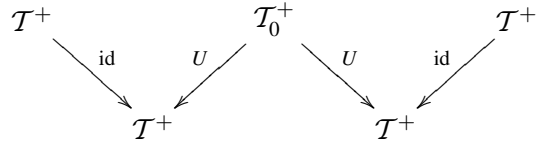
\xymatrix{\calt^{+} \ar[dr]^{\rm id} &&\calt^{+}_{0} \ar[dl]_{U} \ar[dr]^{U}&&
\calt^{+} \ar[dl]_{\rm id}\\
&\calt^{+}& &\calt^{+}}$$ From the absolute grading determination, we see
immediately that
$$d(S^3_{\smash{p/q}}(K), \fraks_0)=d(S^3_{\smash{p/q}}(U),\fraks_0)-2,$$
since the lowest degree term in $\calt^+_0\cong H_*(A^+_0)$ is killed by
multiplication with $U$ and whose degree is 2 lower than that of the element
mapping to the lowest degree element in $\calt^+\cong H_*(B^+)$.

When $q$ is even (see the complex for $S^3_{5/4}(K)$ as an example), the
central subcomplex has the form
$$
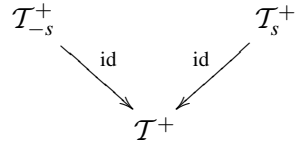
\xymatrix{\calt^{+}_{-s} \ar[dr]^{\rm id}&&
\calt^{+}_s\ar[dl]_{\rm id}\\
&\calt^{+}}$$ where $s=\lfloor(p+q+1)/(2q)\rfloor$. Note that here we need the
condition $q\leq p$ so that the two maps above are identities. So by the
absolute grading determination, we obtain
$$d(S^3_{\smash{p/q}}(K), \fraks_0)=d(S^3_{\smash{p/q}}(U),\fraks_0).\proved$$\end{proof}

{\bf Remark}\qua \fullref{UnknotTrefoil} could also be obtained
by a surgery formula on correction terms given by Rustamov
\cite[Proposition 3.1]{Rustamov}. Note that the spin \spinc structure $\fraks_0$ corresponds to
$l\equiv(p-1)(1-x)/2\pmod{p}$, where $qx\equiv-1\pmod{p}$ in his identification
of the \spinc structures of $S^3_{\smash{p/q}}(K)$ with $\bbz/p\bbz$.
\newpage

\begin{prop}\label{CorrectionTermComputation} For a fake right-handed trefoil
knot $K$, we have the following computation of correction terms:
\begin{align*}
d(S^3_{(18k+9)/(3k+1)}(K), \fraks_0)=&\begin{cases}0, &\text{if $k$ is even},\\
\!+\frac12,&\text{if $k$ is odd.}\end{cases}\\
d(S^3_{(18k+9)/(3k+2)}(K), \fraks_0)=&\begin{cases}0, &\text{if $k$ is even},\\
-\frac12,&\text{if $k$ is odd.}\end{cases}\\
d(S^3_{(18k+9)/(15k+7)}(K), \fraks_0)=&\begin{cases}-2, &\text{if $k$ is even},\\
-\frac32,&\text{if $k$ is odd.}\end{cases}\\
d(S^3_{(18k+9)/(15k+8)}(K), \fraks_0)=&\begin{cases}-2, &\text{if $k$ is even},\\
-\frac52,&\text{if $k$ is odd.}\end{cases}
\end{align*}
Here $\fraks_0$ is the unique spin \spinc structure of the surgered
three-manifold. And further, except for possibly $S^3_{(18k+9)/(3k+1)}(K)\cong
- S^3_{(18k+9)/(3k+2)}(K)$, there are no other homeomorphisms between these
manifolds.
\end{prop}

\begin{proof} To compute the correction terms, we first compute that of the
corresponding surgered manifolds on the unknot $U$. We use notation as in
\fullref{hfhprelim}. For a linearly plumbed four-manifold, let $[S_i]$ be
the homology class of the $i$--th exceptional sphere and $\phi_i=PD([S_i])$.
Since $p=18k+9$ is odd, there is a unique spin \spinc structure, which, by
abuse of notation, we always denote by $\fraks_0$.

The lens space $S^3_{-9}(U)$ is the boundary of the linearly plumbed
four-manifold $X=(-9)$. The \spinc structure $\frakt_0$ on $X$ with
$c_1(\frakt_0)=\phi_0$ restricts to the spin \spinc structure $\fraks_0$ on
$Y$. Hence $\ch(X, \fraks_0)$ contains all cohomology classes of the form
$(2n+1)\phi_0$. We have
$$((2n+1)\phi_0)^2=(2n+1)^2\phi_0^2=-9(2n+1)^2\leq-9,$$
$$d(S^3_{-9}(U),\fraks_0)=\max_{K\in\ch(X,\fraks_0)}\frac{K^2+\rk H^2(X)}4=\frac{(-9)+1}4=-2.\leqno{\hbox{hence}}$$
The lens space $S^3_{-(18k+9)/(3k+1)}(U)$ ($k\geq1)$ can be thought of the
boundary of the linearly plumbed four-manifold
$$X=(-7, -2, \cdots, -2,-4),$$
\newpage
where the number of $-2$ spheres is $k-1$. If $k$ is even, $\phi_0+ \phi_2+
\cdots+ \phi_k$ is characteristic and
$$\ch(X,\fraks_0)=\bigg\{\sum_{i=0}^{k}C_i\phi_i\ {\bf \Big|}\ C_0, C_2,\cdots
C_k\text{ are odd while $C_1, C_3,\cdots C_{k-1}$ are even.}\bigg\}$$ For
$K=C_0\phi_0+\cdots+C_k\phi_k\in\ch(X,\fraks_0)$, we have
\begin{align*}
K^2=&C_0^2\phi_0^2+2C_0C_1\phi_0\phi_1 + C_1^2\phi_1^2+2C_1C_2\phi_1\phi_2\!+
\cdots+\!2C_{k{-}1}C_k\phi_{k{-}1}\phi_k+C_k^2\phi_k^2\\
=&-7C_0^2 + 2C_0C_1-2C_1^2+2C_1C_2-2C_2^2 +\cdots +2C_{k-1}C_k-4C_k^2\\
=&-6C_0^2 -(C_0-C_1)^2-(C_1-C_2)^2-\cdots-(C_{k-1}-C_k)^2-3C_k^2.
\end{align*}
Now $C_0$, $C_k$ and all $C_i-C_{i+1}$ are odd numbers, so we get
$$K^2\leq-6-1-\cdots-1-3=-9-k,$$
and the equality can be obtained by setting the constants $C_0=C_2=\cdots=C_k=1$ and
$C_1=C_3=\cdots=C_{k-1}=0$. Hence for $k$ even, we get
$$d(S^3_{-(18k+9)/(3k+1)}(U),\fraks_0)=\!\max_{K\!\in \ch(X,\fraks_0)}\!\!\!\!\frac{K^2{+}\rk H^2(X)}4
=\frac{(-9{-}k){+}(k{+}1)}{4}=-2.$$ If $k$ is odd, then $\phi_1+\phi_3+\cdots+\phi_k$
is characteristic and
$$\ch(X,\fraks_0)=\bigg\{\sum_{i=0}^{k}C_i\phi_i\ {\bf \Big|}\ C_0, C_2,\cdots
C_{k-1}\text{ are even while $C_1, C_3,\cdots C_k$ are odd.}\bigg\}$$ As in
the case where $k$ is even, for $K=C_0\phi_0+\cdots+C_k\phi_k\in
\ch(X,\fraks_0)$, we get
\begin{align*}
K^2=&-6C_0^2 -(C_0-C_1)^2-(C_1-C_2)^2-\cdots-(C_{k-1}-C_k)^2-3C_k^2 \\
\leq&-1\cdots-1-3=-k-3.
\end{align*}
Note that here we have $C_0$ is even and $C_k, C_i-C_{i+1}$ are odd, and the
equality is obtained by letting $C_0=C_2=\cdots C_{k-1}=0$ and $C_1=C_3=\cdots=
C_k=1$. And the correction term
$$d(S^3_{-(18k{+}9)/(3k{+}1)}(U),\fraks_0)=\!\max_{K\!\in\ch(X,\fraks_0)}\!\!\!\!\frac{K^2{+}\rk H^2(X)}4
=\frac{(-3{-}k){+}(k{+}1)}{4}=-\frac12.$$ Altogether, we have
$$d(S^3_{-(18k+9)/(3k+1)}(U), \fraks_0)=\begin{cases}-2, &\text{if $k$ is even},\\
-\frac12,&\text{if $k$ is odd.}\end{cases}$$
\newpage
By the properties of the
corrections terms \eqref{CorrectionTermProperties}, we get
$$d(S^3_{(18k+9)/(3k+1)}(U), \fraks_0)=\begin{cases}2, &\text{if $k$ is even},\\
\frac12,&\text{if $k$ is odd.}\end{cases}$$
The correction term $d(S^3_{\smash{(18k+9)/(3k+2)}}(U),\fraks_0)$ could be computed by
thinking of the lens space $S^{\smash{3}}_{\smash{-(18k+9)/(3k+2)}}(U)$ as the boundary of
$X=(-5,-2)$ if $k=0$ and of $X=(-6, -k-1, -3)$ if $k>0$. We get
$$d(S^3_{(18k+9)/(3k+2)}(U), \fraks_0)=\begin{cases}0, &\text{if $k$ is even},\\
\frac32,&\text{if $k$ is odd.}\end{cases}$$
The correction term $d(S^3_{\smash{(18k+9)/(15k+7)}}(U),\fraks_0)$ is computed by
thinking of the lens space $S^3_{-(18k+9)/(15k+7)}(U)$ as the boundary of
$$X=\begin{cases}(-2,-2,-2,-3), &\text{if $k=0$,}\\
(-2,-2,-2,-2,-4,-2), &\text{if $k=1$,} \\
(-2,-2,-2,-2,-3,-3,-2), &\text{if $k=2$,}\\
(-2,-2,-2,-2,-3,-2,\cdots,-2,-3,-2), &\text{if $k\geq3$,}
\end{cases}$$
where when $k\geq3$, the number of $(-2)$'s between each $-3$ is $k-2$. We
get
$$d(S^3_{(18k+9)/(15K+7)}(U), \fraks_0)=\begin{cases}0, &\text{if $k$ is even},\\
-\frac32,&\text{if $k$ is odd.}\end{cases}$$
The correction term $d(S^3_{\smash{(18k+9)/(15k+8)}}(U),\fraks_0)$ is computed by
viewing the lens space $S^3_{-(18k+9)/(15k+8)}(U)$ as the boundary of
\begin{gather*}X=(-2,-2,-2,-2,-2,-(k+2),-2,-2).\\ \tag*{\hbox{We have}}
d(S^3_{(18k+9)/(15K+8)}(U), \fraks_0)=\begin{cases}-2, &\text{if $k$ is even},\\
-\frac12,&\text{if $k$ is odd.}\end{cases}\end{gather*}
Now the correction term for $S^3_{\smash{p/q}}(K)$ follows from \fullref{UnknotTrefoil}.

By the properties of the correction terms
\eqref{CorrectionTermProperties}, the only possible homeomorphisms are
\begin{align*}
S^3_{(18k+9)/(3k+1)}(K)\cong&-S^3_{(18k+9)/(3k+2)}(K), &\text{for $k\geq0$,\ \,}\\
S^3_{(18k+9)/(3k+1)}(K)\cong&S^3_{(18k+9)/(3k+2)}(K), &\text{for $k$ even,}\\
S^3_{(18k+9)/(15k+7)}(K)\cong&S^3_{(15k+8)/(15k+8)}(K), &\text{for $k$ even.}
\end{align*}
The last two possibilities can easily be ruled out by their Casson--Walker
invariants.
\end{proof}

Our proof of \fullref{genusonecosmetic} will use the Reidemeister torsion
of a three-manifold. Refer to Milnor \cite{Milnor} or Turaev \cite{Turaev2} for introductions
on torsions. Turaev \cite{Turaev1} and independently Sakai \cite{Sakai}
computed the Reidemeister torsion of $S^3_{\smash{p/q}}(K$).
The corresponding Reidemeister torsion of $S^3_{\smash{p/q}}(K)$ for any $p$--th root of
unity $\xi\neq1$  is
given by
\begin{equation}\label{reidtorsionformula}
\tau(S^3_{\smash{p/q}}(K),\xi)=\Delta_K(\xi)\cdot(\xi-1)^{-1}(\xi^a-1)^{-1},
\end{equation}
where $a$ is defined by $qa\equiv 1\pmod{p}$ and $\Delta_K(t)$ is the Alexander
polynomial of $K$. In fact, this result is true for a knot in any homology
sphere.

\begin{proof}[Proof of \fullref{genusonecosmetic}]
Let $K$ be a genus one knot. Suppose $S^3_{\smash{p/q}}(K)\cong\pm
S^3_{\smash{p/q^\prime}}(K)$, where $q\neq q^\prime$. By \fullref{hfkofcosmeticknot}, we have
$$\Delta_K(t)= t -1 + t^{-1}.$$
By a theorem of Turaev \cite[Theorem 9.1]{Turaev2}, for any $p$--th root of unity
$\xi\neq1$, we have (for some $d\in\bbz/p\bbz$)
$$\tau(S^3_{\smash{p/q}}(K), \xi)=\tau(S^3_{\smash{p/q^\prime}}(K),\xi^d).$$
Here the equality is up to multiplication of $p$--th root of unity and sign. We
have $(d,p)=1$, and by \eqref{reidtorsionformula}, we then get
\begin{equation}
\babs{\Big}{\frac{(\xi-1+\xi^{-1})}{(\xi-1)(\xi^a-1)}}=\babs{\Big}{
\frac{(\xi^d-1+\xi^{-d})}{(\xi^d-1)(\xi^{a^\prime d}-1)}}
\end{equation}
for any $p$--th root of unity $\zeta\neq1$, where $a, a^\prime$ are defined by
$qa\equiv1,q^\prime a^\prime\equiv1\pmod{p}$.

By \fullref{numbertheorem}, we have $a\equiv\pm a^\prime\pmod{p}$,
\begin{align*}p&=12&\text{and}\quad&  \{\pm a, \pm a^\prime\}=\{\pm 1, \pm 5\}\pmod{12},\\
p&=18k+9&\text{and}\quad& \{\pm a,\pm a^\prime\}=\{\pm(6k+1),\pm(6k+5)\}\pmod{p}\tag*{\hbox{or}}\end{align*}
for some nonnegative integer $k$. Equivalently, we have $q\equiv\pm q^\prime
\pmod{p}$,
\begin{align*}p&=12&\text{and}&\quad \{\pm q, \pm q^\prime\}=\{\pm1, \pm5\}\pmod{12},\\
p&=18k+9&\text{and}&\quad\{\pm q,\pm q^\prime\}=\{\pm(3k+1),\pm(3k+2)\}\pmod{p}.\tag*{\hbox{or}}\end{align*}
The $p=12$ case can be ruled out by direct computation of Casson--Walker
invariants and \fullref{NoSmallLSpaceSurgery}.

By \fullref{hfkofcosmeticknot}, \fullref{NoNegativeOnRHT} and
\fullref{NoSmallLSpaceSurgery}, for $p=18k+9$, we only need to consider the
surgeries on fake right-handed trefoil knots with $q$ taking the four values:
$3k+1$, $3k+2$, $15k+7$ and $15k+8$, and the surgeries on fake left-handed
trefoil knots with $q$ taking the four values: $-(3k+1)$, $-(3k+2)$, $-(15k+7)$
and $-(15k+8)$.
\newpage

Note that the $(p/q)$--surgery on a fake left-handed trefoil knot $K$ could be
thought of as the $(-p/q)$--surgery on the mirror of $K$, which is a fake right-handed trefoil knot. Now the theorem follows from \fullref{CorrectionTermComputation}.
\end{proof}

{\bf Remark}\qua Note that in the genus one case, we rule out cosmetic
surgeries of slopes with opposite signs by \fullref{hfkofcosmeticknot}.
However, for the genus two knot $K=9_{44}$, $S^3_{\smash{+1}}(K)$ and $S^3_{-1}(K)$
have the same Heegaard Floer homology, though the two surgered manifolds are
distinguished by their hyperbolic volumes \cite{OSzRationalSurgeries}.
Thus, \fullref{hfkofcosmeticknot} can not be established for genus two
knots (just using Heegaard Floer theory). So for genus two knots, our method might be applied to rule out truly cosmetic surgeries of slopes with the same sign, but there still might be truly cosmetic surgeries of slopes with opposite signs.

For knots in higher Seifert genus, our method might not work (even to rule out
truly cosmetic surgeries with same sign slopes). Our method is to write the
Reidemeister torsions of surgered manifolds on a knot admitting $L$--space
surgeries as products of cyclotomic numbers and use number theoretic results
to get restrictions on possible cosmetic surgery slopes (in favor of \fullref{hfkofcosmeticknot}). As is pointed out in Kadokami \cite{Kadokami}, the
$(-2, 3, 7)$--pretzel knot $K$ has Alexander polynomial
$$\Delta_K(t)=t^{10}-t^9+t^7-t^6+t^5-t^4+t^3-t+1,$$
which is not a product of cyclotomic polynomials. And Fintushel and Stern
\cite{FinStern} showed that 18 and 19 surgeries along $K$ yield lens spaces. So we might need more machinery.

\bibliographystyle{gtart}
\bibliography{link}

\begin{thebibliography}{}
\providecommand\bibmarginpar{\leavevmode\marginpar}
\def\urlstyle#1{{\tt #1}}

\bibitem{Apostol}
\textbf{T\,M Apostol}, \emph{Introduction to analytic number theory}, Springer,
  New York (1976) \xox{MR}{0434929}

\bibitem{BleilerHodgsonWeeks}
\textbf{S\,A Bleiler}, \textbf{C\,D Hodgson}, \textbf{J\,R Weeks},
  \href{http://dx.doi.org/10.2140/gtm.1999.2.23} {\emph{Cosmetic surgery on
  knots}}, from: ``Proceedings of the Kirbyfest (Berkeley, CA, 1998)'', Geom.
  Topol. Monogr. 2, Geom. Topol. Publ., Coventry (1999)  23--34
  \xox{MR}{1734400}

\bibitem{Dvornicich}
\textbf{R Dvornicich}, \emph{On an equation in cyclotomic numbers}, Acta Arith.
  98 (2001) 71--94 \xox{MR}{1831457}

\bibitem{Eftekhary}
\textbf{E Eftekhary}, \href{http://dx.doi.org/10.2140/agt.2005.5.1389}
  {\emph{Longitude {F}loer homology and the {W}hitehead double}}, Algebr. Geom.
  Topol. 5 (2005) 1389--1418 \xox{MR}{2171814}

\bibitem{Ennola}
\textbf{V Ennola}, \href{http://dx.doi.org/10.1016/0022-314X(72)90050-9}
  {\emph{On relations between cyclotomic units}}, J. Number Theory 4 (1972)
  236--247 \xox{MR}{0299585}

\bibitem{FinStern}
\textbf{R Fintushel}, \textbf{R\,J Stern},
  \href{http://dx.doi.org/10.1007/BF01161380} {\emph{Constructing lens spaces
  by surgery on knots}}, Math. Z. 175 (1980) 33--51 \xox{MR}{595630}

\bibitem{Franz}
\textbf{W Franz}, \emph{Uber die Torsion einer \"Uberdeckung}, J. Reine Angew.
  Math. 173 (1935) 245--254

\bibitem{Gabai}
\textbf{D Gabai}, \emph{Foliations and the topology of 3--manifolds II}, J.
  Differential Geom. 26 (1987) 461--478 \xox{MR}{910017}

\bibitem{Ghiggini}
\textbf{P Ghiggini}, \emph{Knot Floer homology detects genus-one fibred links}
  \xox{arXiv}{math.GT/0603445}

\bibitem{GordonLuecke}
\textbf{C\,M Gordon}, \textbf{J Luecke}, \emph{Knots are determined by their
  complements}, J. Amer. Math. Soc. 2 (1989) 371--415 \xox{MR}{965210}

\bibitem{Hedden}
\textbf{M Hedden}, \href{http://dx.doi.org/10.2140/agt.2005.5.1197} {\emph{On
  knot {F}loer homology and cabling}}, Algebr. Geom. Topol. 5 (2005) 1197--1222
  \xox{MR}{2171808}

\bibitem{Neuman}
\textbf{M Jankins}, \textbf{W\,D Neumann}, \emph{Lectures on {S}eifert
  manifolds}, Brandeis Lecture Notes 2, Brandeis University, Waltham, MA (1983)
  \xox{MR}{741334}

\bibitem{Kadokami}
\textbf{T Kadokami}, \emph{Reidemeister torsion of homology lens spaces}, from:
  ``Proceedings of the east asian school of knots, links, and related topics
  (Seoul, Korea)'' (2004)

\bibitem{Kirby}
\textbf{R Kirby}, \emph{Problems in low-dimensional topology}, from:
  ``Geometric topology (Athens, GA, 1993)'', AMS/IP Stud. Adv. Math. 2, Amer.
  Math. Soc., Providence, RI (1997)  35--473 \xox{MR}{1470751}

\bibitem{KMOSz}
\textbf{P\,B Kronheimer}, \textbf{T Mrowka}, \textbf{P\,S Ozsv{\'a}th},
  \textbf{Z Szab{\'o}}, \emph{Monopoles and lens space surgeries}
  \xox{arXiv}{math.GT/0310164}

\bibitem{Lackenby}
\textbf{M Lackenby}, \href{http://dx.doi.org/10.1090/S0894-0347-97-00241-5}
  {\emph{Dehn surgery on knots in 3--manifolds}}, J. Amer. Math. Soc. 10 (1997)
  835--864 \xox{MR}{1443548}

\bibitem{Mathieu}
\textbf{Y Mathieu}, \href{http://dx.doi.org/10.1142/S0218216592000161}
  {\emph{Closed 3--manifolds unchanged by {D}ehn surgery}}, J. Knot Theory
  Ramifications 1 (1992) 279--296 \xox{MR}{1180402}

\bibitem{Milnor}
\textbf{J Milnor}, \emph{Whitehead torsion}, Bull. Amer. Math. Soc. 72 (1966)
  358--426 \xox{MR}{0196736}

\bibitem{Ni3Braids}
\textbf{Y Ni}, \emph{Closed 3--braids are nearly fibred}
  \xox{arXiv}{math.GT/0510243}

\bibitem{Nifibred}
\textbf{Y Ni}, \emph{Knot Floer homology detects fibred knots}
  \xox{arXiv}{math.GT/0607156}

\bibitem{NiSutured}
\textbf{Y Ni}, \href{http://dx.doi.org/10.2140/agt.2006.6.513} {\emph{Sutured
  {H}eegaard diagrams for knots}}, Algebr. Geom. Topol. 6 (2006) 513--537
  \xox{MR}{2220687}

\bibitem{OSzRationalSurgeries}
\textbf{P Ozsv{\'a}th}, \textbf{Z Szab{\'o}}, \emph{Knot {F}loer homology and
  rational surgeries} \xox{arXiv}{math.GT/0504404}

\bibitem{OSzAbsoluteGrade}
\textbf{P Ozsv{\'a}th}, \textbf{Z Szab{\'o}}, \emph{Absolutely graded {F}loer
  homologies and intersection forms for four-manifolds with boundary}, Adv.
  Math. 173 (2003) 179--261 \xox{MR}{1957829}

\bibitem{OSzPlumb}
\textbf{P Ozsv{\'a}th}, \textbf{Z Szab{\'o}},
  \href{http://dx.doi.org/10.2140/gt.2003.7.185} {\emph{On the {F}loer homology
  of plumbed three-manifolds}}, Geom. Topol. 7 (2003) 185--224
  \xox{MR}{1988284}

\bibitem{OSzSurvey}
\textbf{P Ozsv{\'a}th}, \textbf{Z Szab{\'o}}, \emph{Heegaard diagrams and
  holomorphic disks}, from: ``Different faces of geometry'', Int. Math. Ser.
  (N. Y.) 3, Kluwer/Plenum, New York (2004)  301--348 \xox{MR}{2102999}

\bibitem{OSzClosed3M}
\textbf{P Ozsv{\'a}th}, \textbf{Z Szab{\'o}},
  \href{http://projecteuclid.org/getRecord?id=euclid.annm/1105737568}
  {\emph{Holomorphic disks and topological invariants for closed
  three-manifolds}}, Ann. of Math. $(2)$ 159 (2004) 1027--1158
  \xox{MR}{2113019}

\bibitem{OSzMutation}
\textbf{P Ozsv{\'a}th}, \textbf{Z Szab{\'o}},
  \href{http://dx.doi.org/10.1016/j.topol.2003.09.009} {\emph{Knot {F}loer
  homology, genus bounds, and mutation}}, Topology Appl. 141 (2004) 59--85
  \xox{MR}{2058681}

\bibitem{deRham}
\textbf{G de~Rham}, \textbf{S Maumary}, \textbf{M\,A Kervaire}, \emph{Torsion
  et type simple d'homotopie}, Lecture Notes in Mathematics 48, Springer,
  Berlin (1967) \xox{MR}{0222893}

\bibitem{Rong}
\textbf{Y\,W Rong}, \emph{Some knots not determined by their complements},
  from: ``Quantum topology'', Ser. Knots Everything 3, World Sci. Publ., River
  Edge, NJ (1993)  339--353 \xox{MR}{1273583}

\bibitem{Rustamov}
\textbf{R Rustamov}, \emph{The renormalized Euler characteristic and L-space
  surgeries} \xox{arXiv}{math.GT/0506320}

\bibitem{Sakai}
\textbf{T Sakai}, \emph{Reidemeister torsion of a homology lens space}, Kobe J.
  Math. 1 (1984) 47--50 \xox{MR}{784347}

\bibitem{Turaev1}
\textbf{V\,G Turaev}, \emph{Reidemeister torsion and the {A}lexander
  polynomial}, Mat. Sb. $($N.S.$)$ 18(66) (1976) 252--270 \xox{MR}{0433462}

\bibitem{Turaev2}
\textbf{V Turaev}, \emph{Introduction to combinatorial torsions}, Lectures in
  Mathematics ETH Z\"urich, Birkh\"auser Verlag, Basel (2001) \xox{MR}{1809561}

\bibitem{Washington}
\textbf{L\,C Washington}, \emph{Introduction to cyclotomic fields}, second
  edition, Graduate Texts in Mathematics 83, Springer, New York (1997)
  \xox{MR}{1421575}

\end{thebibliography}

\end{document}